\documentclass{article}
\usepackage[T2A]{fontenc}

\usepackage[english]{babel}
\usepackage[tbtags]{amsmath}
\usepackage{amsfonts,amssymb,mathrsfs,amscd}
\usepackage{amssymb, amsmath, dsfont, srcltx}

\newtheorem{theorem}{Theorem}[section]

\newtheorem{e-proposition}[theorem]{Proposition}
\newtheorem{corollary}[theorem]{Corollary}
\newtheorem{e-definition}[theorem]{Definition\rm}

\def\a{\alpha}
\def\b{\beta}
\def\om{\omega}
\def\e{\varepsilon}
\def\g{\gamma}
\def\G{\Gamma}
\def\l{\lambda}
\def\p{\partial}
\def\D{\Delta}
\def\Om{\Omega}

\def\vr{\varrho}

\def\L{\Lambda}

\def\H1{W^{1,2}}
\def\H2{W^{2,2}}

\def\di{\,d}

\newcommand{\EE}{\mathbb{E}}
\newcommand{\NN}{\mathds{N}}

\newcommand{\PP}{\mathbb{P}}

\newcommand{\RR}{\mathds{R}}

\def\cL{\mathcal{L}}
\def\Op{\mathcal{H}}

\def\bc{\mathcal{B}}

\def\spec{\sigma}

\renewcommand\e{\varepsilon}

\newcommand{\supp}{\operatorname{supp}}
\newcommand{\dist}{\operatorname{dist}}

\newcommand{\randfunct}{f}

\renewcommand{\leq}{\leqslant}
\renewcommand{\geq}{\geqslant}

\begin{document}

\title{Spectral localization for quantum Hamiltonians with weak random delta interaction}

\author{Denis I. Borisov$^1$, Matthias T\"aufer$^2$, Ivan Veseli\'c$^2$}

\date{\empty}

\maketitle
{\small 
\begin{quote}
1) Department of Differential Equations, Institute of Mathematics with Computer Center, Ufa Federal Research Center, Russian Academy of Sciences, Chernyshevsky. st.~112, Ufa, 450008,
Russia
\\
Faculty of Physics and Mathematics, Bashkir
State Pedagogical University, October rev. st.~3a, Ufa, 450000,
Russia 
\\ 
Faculty of Natural Sciences, University of Hradec Kr\'alov\'e, Rokitansk\'eho 62, 500 03, Hradec Kr\'alov\'e, Czech Republic
\\
2) Fakult\"at f\"ur Mathematik, Technische Universit\"at Dortmund, 44227 Dortmund, Germany
\\
emails: borisovdi@yandex.ru, m.taeufer@qmul.ac.uk, iveselic@math.tu-dortmund.de
\end{quote}
 }

\begin{abstract}

We consider a negative Laplacian in multi-dimensional Euclidean space (or a multi-dimensional layer)
with a weak disorder random perturbation.
The perturbation consists of a sum of lattice translates of a
delta interaction supported on a compact manifold of co-dimension one
and modulated by coupling constants which are
independent identically distributed random variables times a small disorder parameter.
We establish that the spectrum of the considered operator is almost surely a fixed set, characerize its minimum,
give an initial length scale estimate and the Wegner estimate
and conclude that there is a small zone of a pure point spectrum containing the almost sure spectral bottom.
The length of this zone is proportional to the small disorder parameter. 
\end{abstract}

\section{Introduction}

One of the most prominent aspects in the theory of random Hamiltonians is Anderson localization.
From the physical point of view, localization in a disordered medium describes the situation when waves within a certain energy or frequency zone do not propagate through the medium.
Mathematically, this means in particular that the spectrum of the considered random operator contains
a region where almost surely continuous spectrum is absent and consequently the entire spectrum
consists of the closure of the set of eigenvalues.

In this paper, we give the key results and ingredients for Anderson localization in a new model, namely, random delta-interactions supported on surfaces or, more precisely, on manifolds of co-dimension one.
Previously, Hamiltonians with random delta-interactions have been studied only in the case where the support of the singular interactions is discrete, leading to Hamiltonians with point interactions (in low dimensions).
To the best of our knowledge the only two results valid in dimensions larger than one are \cite{HislopKK-05}, where Hislop, Kirsch, and Krishna establish localization for random point interactions in dimension up to $3$, and \cite{KloppP-08} where Klopp and Pankrashkin study quantum graphs with random vertex couplings.

The study of localization for Hamiltonians with random singular interactions poses several challenges which are not present for simple models like the Anderson or alloy type model.
As the very name says, the perturbation is no longer a potential but singular, in fact, a measure.
Furthermore, at least in the generality which we are treating, it is not required that the perturbation be monotone with respect to the coupling constants.
It is well known that non-monotonicity in the random variables requires adapted methods in the context of localization.
One possibility to 'tame' the singular interaction is to transform the operator in to another one which is spectrally equivalent, but easier to analyze.
This is in fact the first step of our proof, however one has to pay a price for this.
The transformed operator does no longer depend in a linear fashion on the coupling constants, but in a non-linear one.

The model which we present here is just one instance of a general class of models which we are able to treat and which we will discuss in detail in a subsequent paper \cite{BTV}.
The first steps to implement the multiscale analysis proof of localization for a general random operator was initiated in the works \cite{BGV-16}, \cite{B1}.
In these papers, the initial length scale estimate at the bottom of the spectrum was proved for the Schr\"odinger operator in a multi-dimensional layer with a small random perturbation distributed in the cells of some periodic lattice.
The perturbation was described by an abstract operator of the form $\cL(t)=t \cL_1+t^2 \cL_2 +t^3 \cL_3(t)$ depending on a real parameter $t$ with the coefficients $\cL_i$ satisfying certain rather weak conditions.
This operator was considered at each periodicity cell and the parameter $t$ was replaced then by $\e\om_k$, where $\e$ is a global small parameter and $\om_k$ is a random variable associated with the cell.
By choosing various particular cases of the operators $\cL_i$, the proposed model covered many interesting particular examples, both known and new.

A complementary ingredient to complete the multiscale analysis and conclude localization is a Wegner estimate.
Since we are working in the weak disorder regime, both the initial length scale and the Wegner estimate we aim for are restricted to certain zones is the energy $\times$ disorder diagram.
For this reason it is necessary to resolve a preliminary question: Is there any spectrum in the zone where our Wegner estimate applies? Of course it is trivial to prove a Wegner estimate for energies in the resolvent set, but this cannot lead to a localization result in the proper sense of the word.

So the question is, whether the spectrum expands under the influence of the random perturbation or not. And if yes: At which rate?
This problem was successfully solved in \cite{BVECCM18} for the above described general model.
This enables us to turn to the final missing step, namely, the Wegner estimate.
By employing and  modifying the technique proposed in \cite{HislopK-02}, we succeeded to reprove the Wegner estimate.
Depending on the rate of the expansion of the spectrum in the weak disorder regime the Wegner estimate takes on different forms.
The results for the general model will be presented and proved in \cite{BTV}.

As an excerpt and illustration of the aforementioned general results we treat
in the  paper a new interesting example, namely, a weak random delta interaction in
multi-dimensional Euclidean space (or a layer) where the delta-interaction is supported on a manifold of codimension one.
We present an initial length scale and a Wegner estimate for this model and as a corollary, Anderson localization at the bottom of the spectrum.

\section{Problem and results}
Our study applies both to operators defined on the whole Euclidean space as well as to ones on strips, layers, or higher-dimensional analogs.
In fact, the results for operators on the whole Euclidean space $\RR^n $ can be recovered as a corollary
from the ones for operators on layer domains in $\RR^{n+1}$ using the reduction described in \cite[\S 3.7]{BGV-16}.
For this reason we discuss in the present note strip domains only.
 Let $x'=(x_1,\ldots,x_n)$, $x=(x',x_{n+1})$, $n\geqslant 1$, be Cartesian coordinates in $\RR^n$ and $\RR^{n+1}$, respectively, $\Pi$ be the multidimensional layer $\Pi:=\{x:\, 0<x_{n+1}<d\}$
of width $d>0$. In $\RR^n$ we introduce a periodic lattice $\G$ with a basis $e_1$, \ldots, $e_n$ and the unit cell
of this lattice is denoted
by $\square'$, i.e., $\square':=\{x': x'=\sum\limits_{i=1}^{n} t_i e_i,\; t_i\in(0,1)\}$. We also define $\square:=\square'\times(0,d)$.

Let $\om:=\{\om_k\}_{k\in\G}$ be a sequence of independent identically distributed random variables with distribution measure $\mu$ with an absolutely continuous probability density $h$.
We assume that the support of $\mu$ satisfies
\begin{equation}
 \label{eq:support_mu}
 -1\leqslant a=\min\supp\mu<\max\supp\mu=1
\end{equation}
and we introduce the product probability measure $\PP=\bigotimes_{k\in\G} \mu$ on the configuration space $\Om:=\times_{k\in\G} [-1,1]$;
the elements of this space are sequences $\om:=\{\om_k\}_{k\in\G}$. By $\EE(\cdot)$ we denote the expectation value of a random variable w.r.t.~$\PP$.
Let us comment on our assumption that the probability density function $h$ is absolutely continuous.
All proofs of Wegner estimates for non-monotone random models in multidimensional continuum space known to us require some form of differentialbility
of the density $h$.
This is for instance the case for alloy-type potentials with sign-changing single site potentials, see e.~g.~\cite{HislopK-02} and \cite{Veselic-02a},
and for the random displacement model, see e.~g.~\cite{KloppLNS-12}. Our proof of the Wegner estimate uses the ideas of~\cite{HislopK-02}, and thus inherits its regularity requirement.

Let $M_0$ be a closed bounded $C^3$ manifold in $\square$ of codimension one. The outward normal vector to $M_0$ is denoted by $\nu$. The manifold $M_0$ is assumed to be separated
from the boundary $\p\square$ by a positive distance. By $M_k$, $k\in\G$, we denote the translation of $M_0$ along $\G$:
\begin{equation*}
M_k:=\{x\in\Pi:\, x-k\in M_0\}
\end{equation*}
and by $M:=\bigcup\limits_{k\in\G} M_k$ their union.

By $y=(y_1,\ldots,y_{n-1})$ we denote some local coordinates on $M_0$, while $\vr$ is the distance from a point in $\Pi$ to $M_0$ measured along $\nu$.
We will need this parametrization only locally.
Since $M_k$ are translations of $M_0$, the coordinates $(y,\vr)$ are well-defined in a neighborhood of each $M_k$, $k\in\G$,
and hence, in a neighborhood of the entire set $M$. Let $\randfunct=\randfunct(y,t)$ be a real function on $M_0\times[-t_0,t_0]$ for some $t_0>0$.
We assume that
 $\randfunct\in C^4(M_0\times[-t_0,t_0])$. We extend the function $\randfunct$ periodically to the entire set $M$.

Our random operator is introduced as the  Laplacian with random  delta interaction with  a small coupling parameter.
This is the operator in $L_2(\Pi)$ with the differential expression
\begin{equation}\label{3.10}
\Op^\e(\om):=-\D
\quad\text{in}\quad L^2(\Pi).
\end{equation}
The domain of the operator $\Op^\e(\om)$  consists of the functions $u\in \H2(\Pi\setminus M)$ satisfying the boundary conditions
\begin{equation}\label{3.11}
[u]_{M_k}=0,\quad \left[\frac{\p u}{\p\varrho}\right]_{M_k}=-\e \om_k \randfunct(\cdot,\e\om_k) u\big|_{M_k},\quad [u]_{M_k}:=u\Big|_{\genfrac{}{}{0 pt}{}{\vr=+0}{y\in M_k}}-u\big|_{\genfrac{}{}{0 pt}{}{\vr=-0}{y\in M_k}},
\end{equation}
where $\e$ is a (sufficiently small) positive parameter.
The most important case is when $f$ does not depend on the second argument.
Then, the second condition in~\eqref{3.11} reduces to
\[
 \left[\frac{\p u}{\p\varrho}\right]_{M_k}=-\e \om_k \randfunct(\cdot) u\big|_{M_k}
\]
and the single site perturbation depends linearly on the coupling $\e \om_k$.
On $\p\Pi$ we impose either the  Dirichlet or Neumann condition.
We write the latter condition as
\begin{equation}\label{eq:boundary-condition-B}
\bc u=0
\end{equation}
on $\p\Pi$, and $\bc u=u$ or $\bc u=\frac{\p u}{\p x_{n+1}}$.
We stress that on the upper and lower boundaries of $\p\Pi$ we can have different boundary conditions.

The function
$\randfunct(y,t)$  can change its sign and there is no assumption on the type of dependence on $t$ apart from the required regularity. In particular, the considered random operator can be non-monotone and non-linear with respect to the random variables.

To formulate our main results, we need to introduce additional notation. First of all, for $\eta \in [- t_0, t_0]$ we define the auxiliary operator
\begin{equation}\label{2.4}
\Op^\eta_\square:=-\Delta
\end{equation}
in $L^2(\square)$ on the domain consisting of functions $u$ in the Sobolev space $\H2(\square\setminus M_0)$, which satisfy condition (\ref{eq:boundary-condition-B}) on $\p\square\cap\p\Pi$,
condition (\ref{3.11}) with $\e \omega_0 = \eta$ on $M_0$
 and periodic boundary conditions on  $\g:=\p\square\setminus\p\Pi$.
The operator $\Op^\eta_\square$ is self-adjoint.

By $\L^\eta$ we denote the lowest eigenvalue of the operator $\Op^\eta_\square$ and by $\Psi^\eta$ the associated positive eigenfunction.
The symbol $\L_0$ stands for the smallest eigenvalue of the operator
\begin{equation*}
-\frac{d^2}{dx_{n+1}^2}
\quad \text{on}\quad (0,d)
\end{equation*}
subject to boundary condition (\ref{eq:boundary-condition-B}).  The associated eigenfunction is denoted by
$\Psi_0=\Psi_0(x_{n+1})$ and it is normalized as follows:
\begin{equation*}
\|\Psi_0\|_{L^2(0,d)}=\frac{1}{\sqrt{|\square'|}},
\end{equation*}
where $|\square'|$ stands for the volume of $\square'$.
Denote
\begin{equation*}
\L_1=\int\limits_{M_0} \randfunct(y,0)\Psi_0^2(x(y,0))\di y.
\end{equation*}
Our \emph{main assumption} is
\begin{equation}\label{1.0}
\L_1\ne 0.
\end{equation}
This means that while the random perturbations are non-monotone,
there is a dominating sign with respect to the unperturbed ground state.
Our first main result describes the location of the spectrum of $\Op^\e(\om)$,
which will be denoted by $\spec(\Op^\e(\om))$.

\begin{theorem}\label{th2.1}
For all sufficiently small $\e$
there exists a closed set $\Sigma_\e\subset\RR$ such that
\begin{equation}\label{2.18}
\spec(\Op^\e(\om))=\Sigma_\e\quad\PP-a.s.
\end{equation}
The set $\Sigma_\e$ is equal to the closure of the union of spectra for all periodic realizations of $\Op^\e(\cdot)$:
\begin{equation}\label{2.19}
\Sigma_\e=\overline{\bigcup\limits_{N\in \NN} \quad \bigcup\limits_{\xi \ \text{is $2^N \G$-periodic}} \spec\big(\Op^\e(\xi)\big)},
\end{equation}
where the second union is taken over all sequences $\xi: \G\to\supp \mu $, which are periodic  with respect to   the sublattice $2^N \G:=\{2^N q:\, q\in\G\}$.

The  minimum of $\Sigma_\e$ is given by the identity
\begin{equation}\label{2.20}
\min\Sigma_\e=\L^{\e_*},
\end{equation}
where
\begin{equation*}
\e_*=\e \quad\text{if}\quad \L_1<0\qquad\text{and}\qquad \e_*=\e a\quad\text{if}\quad \L_1>0.
\end{equation*}
\end{theorem}
Recall that $a = \inf \supp \mu$, cf. Ineq.~\eqref{eq:support_mu}.

Our second main result is the initial length scale estimate at the bottom of the spectrum. First we need some additional notation.
Given $\a\in\G$, $N\in\NN$, we denote
\begin{align*}
&\Pi_{\a,N}:=\bigg\{x:\, x'=\a+\sum\limits_{j=1}^{n} a_j e_j,\ a_j\in(0,N),\ 0<x_{n+1}<d\bigg\},
\\
&\G_{\a,N}:=\bigg\{x'\in\G:\, x'=\a+\sum\limits_{j=1}^{n} a_j e_j,\ a_j=0,1,\ldots,N-1\bigg\}.
\end{align*}

In $L_2(\Pi_{\a,N})$ we introduce the operator
\begin{equation*}
\Op_{N}^\e(\om):=-\D 
\end{equation*}
subject to conditions (\ref{3.11}) on $M\cap\Pi_{\a,N}$, to condition (\ref{eq:boundary-condition-B}) and to the boundary condition
\begin{equation*}
\frac{\p u}{\p\nu}= \rho^{\e_*}u\quad\text{on}\quad \g_{N}:=\p\Pi_{N}\setminus\p\Pi.
\end{equation*}
Here $\nu$ denotes the outward normal to $\g_N$ and
\begin{equation*}
\rho^{\e_*}:=\frac{1}{\Psi^{\e_*}} \frac{\p\Psi^{\e_*}}{\p\nu}.
\end{equation*}
We can prove that the latter function is well-defined, $\G$-periodic and belongs to $L_\infty(\g_N)$.
Denote by $\chi_S$ the characteristic function of a set $S$.
We then have the following initial length scale estimate.
It estimates the probability of the rare event that the lowest eigenvalues of a segment Hamiltonian lies very close to $\min\Sigma_\e$.
Actually, we formulate a stronger statement which is implied by standard Combes-Thomas estimates.

\begin{theorem}
 \label{th4.2}
Let $\tau \in \NN$, $\tau \geq 5$.
There are $c_0, c_1, c_2, N_1 > 0$ such that for all $N \geq N_1$, the set
\begin{equation*}
J_N:=\left[\frac{8}{\sqrt{|\L_1|\mathbb{E}(|\om_0|)}} \frac{1}{N^{\frac{1}{2}}}, \frac{c_0}{N^{\frac{2}{\tau}}}\right]
\end{equation*}
is non-empty and for all $\e \in J_N$, $\a \in \Gamma$, $N > N_1$ as well as all for all $\b_1, \b_2 \in \Gamma_{\a,N}$, and $m_1, m_2 > 0$ such that $B_k:=\Pi_{\b_k,m_k}\subset\Pi_{\a,N}$ for $k \in \{1,2\}$, we have
 \begin{align*}
 \mathbb{P}&\left(\forall \l\leqslant\L^{\e_*}+\frac{1}{2}\left(\frac{\e}{c_0}\right)^{\frac{\tau}{4}}: \, \|\chi_{B_1}(\Op_{\a,N}^\e(\om)-\l)^{-1}\chi_{B_2}\|\leqslant 2
\left(\frac{\e}{c_0}\right)^{-\frac{\tau}{4}} e^{-c_2 \dist(B_1,B_2)
\left(\frac{\e}{c_0}\right)^{\frac{\tau}{4}}}
 \right)
\\
&\geqslant 1-\left(\frac{\e}{c_0}\right)^{-\frac{n(\tau-1)}{2}} e^{-c_1\left(\frac{\e}{c_0}\right)^{-\frac{n}{2}}}.
\end{align*}

\end{theorem}

Our next result is the Wegner estimate.
It estimates the probability to find spectrum in a very small energy interval.

\begin{theorem}\label{thm:Wegner}
 There are $C_1 > 0$,
 $C_2 > 0$, and $\e_0 > 0$ such that for all $0 < \e \leq \e_0$, all $E_0 \leq \Lambda_0 - C_2 \e^2$, all $E \in (- \infty, E_0]$, all $\kappa \leq \tfrac{1}{4} |\L_0 -  E_0|$, all $\alpha \in \Gamma$, and all $N \in \NN$, we have
 \begin{align}
 \label{eq:Wegner}
  \PP ( \dist(\sigma(\Op_{\a,N}^\e(\omega)), E) \leq \kappa )
  \leq
  \frac{C_1\Vert h \Vert_{W^{1,1}(\RR)}}{| \L_0 -  E_0 |}
  \kappa
   d \, N^{2n}.
 \end{align}
\end{theorem}
Note that our Wegner estimate applies only to energy intervals $[E-\kappa, E +\kappa] $ with a security distance to the minimum of the unperturbed spectrum $\Lambda_0$ of the order
$\epsilon^2$. Due to the assumption \eqref{1.0} the spectrum expands proportionally to $\epsilon$ when we switch on the random perturbation. Thus, our Wegner estimate indeed applies
to a region with non-empty intersection with $\Sigma_\epsilon$.

The initial length scale estimate in Theorem~\ref{th4.2} and the Wegner estimate in Theorem~\ref{thm:Wegner} imply via multiscale analysis spectral localization  on the bottom of  $\Sigma_\e$. This is formulated in the next

\begin{corollary}\label{th2.2} Assume (\ref{1.0}). There exists an $\e_0>0$ and a constant $C>0$ such that for all $\e\in(0,\e_0)$ the set
$\Sigma_\e\cap[\L^{\e_*},\L^{\e_*}+C\e]$ is non empty and contains no continuous spectrum.
\end{corollary}

This theorem means that at the bottom of the spectrum $\Sigma_\e$, there is a segment of length $O(\e)$ containing only pure point spectrum.
By using the approach of~\cite[Section 3.17]{BGV-16}, the above results can be extended from operators on the layer $\Pi = \RR^n \times (0,d)$ to operators on the entire space $\RR^n$.

The proofs of Theorems~\ref{th2.1} to~\ref{thm:Wegner} take some twenty pages, so we can just briefly discuss them here.
Following~\cite[Sect. 8.5]{Borisov-07}, one sees that the operator $\Op^\e(\xi)$ is unitarily equivalent to an operator $\tilde \Op^\e(\xi)$ which on every elementary cell $\square + k$, $k \in \Gamma$ acts as
\[
 - \Delta + \e \xi_k \cL_1 + \e^2 \xi_k^2 \cL_2 + \e^3 \xi_k^3 \cL_3(\e \xi_k),
\]
 where $\cL_1$, $\cL_2$ are operators that are relatively bounded with respect to $-\Delta$ and $\cL_3(t)$ is uniformly relatively bounded with respect to $-\Delta$.
Since the spectra of $\Op^\e(\xi)$ and $\tilde \Op^\e(\xi)$ coincide, it suffices to study $\tilde \Op^\e(\xi)$.

The most part of Theorem~\ref{th2.1} was proved in \cite{BVECCM18} whereas identity (\ref{2.20}) is a new result. It is proved by introducing a special Mezincescu condition on the lateral boundaries of the periodicity cells in $\Pi$ as in \cite{BVJFA13}, \cite{B1}, analyzing the behavior of $\L^\eta$ for small $\eta$, and applying the minimax principle.
To prove Theorem~\ref{th4.2}, we adapt the approach used in \cite{BGV-16}, \cite{B1}. The Wegner estimate in Theorem~\ref{thm:Wegner} is a new non-trivial ingredient. Here we employed and extended the technique proposed in \cite{HislopK-02} combining it with results of \cite{BVECCM18}.

\section*{Acknowledgments}

The research of D.I.B. is financially supported by Russian Science Foundation (project no. 17-11-01004).
M.T.~and I.V.~were partially supported by the grant VE 253/6-1 of the Deutsche Forschungsgemeinschaft.

\end{document}